\documentclass[reqno]{amsart}
\usepackage{amsmath, pb-diagram}
\usepackage{amssymb}
\usepackage{amscd}

\begin{document}
\title [VERY
CLEAN MATRICES OVER LOCAL RINGS]
 {VERY
CLEAN MATRICES OVER LOCAL RINGS}

\author{H. Chen}
\address{H. Chen, Department of Mathematics, Hangzhou Normal University, Hangzhou,
China} \email{huanyinchen@aliyun.com}

\author{B. Ungor}
\address{Burcu Ungor, Department of Mathematics, Ankara University, 06100 Ankara, Turkey}
\email{bungor@science.ankara.edu.tr}

\author{S. Halicioglu}
\address{Sait Hal\i c\i oglu, Department of Mathematics, Ankara University, 06100 Ankara, Turkey}
\email{halici@science.ankara.edu.tr}
\date{}

\newtheorem {thm}{Theorem}[section]
\newtheorem{lem}[thm]{Lemma}
\newtheorem{prop}[thm]{Proposition}
\newtheorem{cor}[thm]{Corollary}
\newtheorem{df}[thm]{Definition}
\newtheorem{nota}{Notation}
\newtheorem{note}[thm]{Remark}
\newtheorem{ex}[thm]{Example}
\newtheorem{exs}[thm]{Examples}
\newtheorem{rmk}[thm]{Remark}
\newtheorem{quo}[thm]{Question}

\def\bc{\begin{center}}
\def\ec{\end{center}}
\def\no{\noindent}
\def\hang{\hangindent\parindent}
\def\textindent#1{\indent\llap{[#1]\enspace}\ignorespaces}
\def\re{\par\hang\textindent}

\begin{abstract} An element $a\in R$ is very clean
provided that there exists an idempotent $e\in R$ such that
$ae=ea$ and either $a-e$ or $a+e$ is invertible. A ring $R$ is
very clean in case every element in $R$ is very clean. We explore
the necessary and sufficient conditions under which a triangular
$2\times 2$ matrix ring over local rings is very clean. The very
clean $2\times 2$ matrices over commutative local rings are
completely determined. Applications to matrices over power series
are also obtained.\vskip3mm\no{\it\bf Key Words:}\ \ very clean
ring; very clean matrix; local ring.

\vskip3mm \no{\it\bf 2010 Mathematics Subject Classification:}\ \ 15A13, 15B99, 16L99.
\end{abstract}

\maketitle

\bc{\bf 1. INTRODUCTION}\ec \vskip4mm \no A ring $R$ is strongly
clean provided that for any $a\in R$ there exist an idempotent
$e\in R$ and an element $u\in U(R)$ such that $a=e+u$ and $ae=ea$,
where $U(R)$ is the set of all units in $R$. Recently, strong
cleanness has been extensively studied in the literature (cf.
[2-3], [5] and [6-8]). We say that an element $a\in R$ is very
clean provided that there exists an idempotent $e\in R$ such that
$ae=ea$ and either $a-e$ or $a+e$ is invertible. A ring $R$ is
very clean in case every element in $R$ is very clean. Clearly,
strong cleanness implies the very cleanness. But the converse is
not true (see Lemma 2.4). The motivation of this note is to
explore very clean matrices over local rings, which also extend
weak cleanness from commutative rings to noncommutative rings~(cf.
[1]). We will construct a large class of very clean rings which
are not strongly clean. Let $A$ and $B$ be local rings, let $V$ be
an $A$-$B$-bimodule, and let $R=\{\left(
\begin{array}{ll}
a&v\\
0&b
\end{array}
\right)~| a\in A,b\in B,v\in V\}.$ We prove that $R$ is very clean
if and only if $\frac{1}{2}\in A$ and $\frac{1}{2}\in B$; or $R$
is strongly clean. The characterization of the very cleanness of
$2\times 2$ matrices over commutative local rings are completely
determined. Let $R$ be a commutative local ring, and let
$\varphi\in M_2(R)$. We prove that $\varphi\in M_2(R)$ is very
clean if and only if $\frac{1}{2}\in R$; or $\varphi\in M_2(R)$ is
strongly clean. Let $R$ be a weakly bleached local ring. We
further show that $A(x)\in M_2\big(R[[x]])$ if and only if
$A(0)\in M_2(R)$ is very clean.

Throughout, all rings are associative with an identity. If
$\varphi\in M_n(R)$, we use $\chi (\varphi)$ to stand for the
characteristic polynomial $det( tI_n-\varphi)$. $M_n(R)$ and
$T_n(R)$ denote the ring of all $n\times n$ matrices and  the ring
of all $n\times n$ upper triangular matrices over $R$,
respectively.

\vskip15mm \bc{\bf 2. Triangular Matrix Rings}\ec

\vskip4mm \no A ring $R$ is local in case it has only one maximal
right ideal. As is well known, a ring $R$ is local if and only if
$a+b=1$ in $R$ implies that either $a$ or $b$ is invertible. The
purpose of this section is to consider very cleanness for a kind
of triangular matrices.

\vskip4mm \no{\bf Lemma 2.1.}\ \ {\it Let $A$ and $B$ be local
rings, let $V$ be an $A$-$B$-bimodule, and let $$R=\{\left(
\begin{array}{ll}
a&v\\
0&b
\end{array}
\right)~| a\in A,b\in B,v\in V\}.$$ Then the following are
equivalent:}
\begin{enumerate}
\item [(1)]{\it $R$ is very clean.} \vspace{-.5mm}
\item [(2)]{\it If $a\pm 1\in J(A), b\in J(B)$ or $a\in J(A), b\pm 1\in J(B)$, and $v\in V$,
there exists $x\in V$ such that $ax-xb=v$.}
\end{enumerate}\no{\it Proof.}\ \ $(1)\Rightarrow (2)$ If $a\pm 1\in J(A), b\in J(B)$ and $v\in V$,
then $r=:\left(
\begin{array}{cc}
a&-v\\
0&b
\end{array}
\right)\in R$ is very clean. Thus, we can find an idempotent $e_x$
such that $re_x=e_xr$ and either $r+e_x\in U(R)$ or $r-e_x\in
U(R)$, where $e_x=\left(
\begin{array}{cc}
0&x\\
0&1
\end{array}
\right)$. As $re_x=\left(
\begin{array}{cc}
0&ax-v\\
0&b
\end{array}
\right)=\left(
\begin{array}{cc}
0&xb\\
0&b
\end{array}
\right)=e_xr$. Hence $ax-v=xb$, and so $ax-xb=v$.

If $a\in J(A), b\pm 1\in J(B)$ and $v\in V$, then $r=:\left(
\begin{array}{cc}
a&v\\
0&b
\end{array}
\right)\in R$ is very clean. Thus, we can find an idempotent $f_x$
such that $rf_x=f_xr$ and either $r+f_x\in U(R)$ or $r-f_x\in
U(R)$, where $f_x=\left(
\begin{array}{cc}
1&x\\
0&0
\end{array}
\right)$. As $rf_x=\left(
\begin{array}{cc}
a&ax\\
0&0
\end{array}
\right)=\left(
\begin{array}{cc}
a&v+xb\\
0&0
\end{array}
\right)=f_xr$. Hence, $ax=v+xb$, as required.

$(2)\Rightarrow (1)$ Let $r=\left(
\begin{array}{cc}
a&v\\
0&b
\end{array}
\right)\in R$.

$(i)$ $a\in J(A),b\in J(B)$. Then $r-\left(
\begin{array}{cc}
1&0\\
0&1
\end{array}
\right)\in U(R)$, and so $r\in R$ is very clean.

$(ii)$ $a\not\in J(A),b\not\in J(B)$. Then $r\in U(R)$, and so
$r\in R$ is very clean.

$(iii)$ $a\not\in J(A),b\in J(B)$. If $a\pm 1\in J(A)$, then there
exists $x\in V$ such that $ax-xb=-v$. Hence, $r+\left(
\begin{array}{cc}
0&x\\
0&1
\end{array}
\right)\in U(R)$. One easily checks that $r\left(
\begin{array}{cc}
0&x\\
0&1
\end{array}
\right)=\left(
\begin{array}{cc}
0&ax+v\\
0&b
\end{array}
\right)=\left(
\begin{array}{cc}
0&xb\\
0&b
\end{array}
\right)=\left(
\begin{array}{cc}
0&x\\
0&1
\end{array}
\right)r$. Therefore $r\in R$ is very clean.

If $a+1\not\in J(A)$, then $r+\left(
\begin{array}{cc}
1&0\\
0&1
\end{array}
\right)\in U(R)$; hence, $r\in R$ is very clean.

If $a-1\not\in J(A)$, then $r-\left(
\begin{array}{cc}
1&0\\
0&1
\end{array}
\right)\in U(R)$; hence, $r\in R$ is very clean.

$(iv)$ $a\in J(A), b\not\in J(B)$.  If $b\pm 1\in J(B)$, then there
exists $x\in V$ such that $ax-xb=v$. Hence, $r+\left(
\begin{array}{cc}
1&x\\
0&0
\end{array}
\right)\in U(R)$. One easily checks that $r\left(
\begin{array}{cc}
1&x\\
0&0
\end{array}
\right)=\left(
\begin{array}{cc}
a&ax\\
0&0
\end{array}
\right)=\left(
\begin{array}{cc}
a&v+xb\\
0&0
\end{array}
\right) = \left(
\begin{array}{cc}
1&x\\
0&0
\end{array}
\right)r$. Therefore $r\in R$ is very clean.

If $b+1\not\in J(B)$, then $r+\left(
\begin{array}{cc}
1&0\\
0&1
\end{array}
\right)\in U(R)$; hence, $r\in R$ is very clean.

If $b-1\not\in J(B)$, then $r-\left(
\begin{array}{cc}
1&0\\
0&1
\end{array}
\right)\in U(R)$; hence, $r\in R$ is very clean. \hfill$\Box$

\vskip4mm \no{\bf Theorem 2.2.}\ \ {\it Let $A$ and $B$ be local
rings, let $V$ be an $A$-$B$-bimodule, and let $$R=\{\left(
\begin{array}{ll}
a&v\\
0&b
\end{array}
\right)~| a\in A,b\in B,v\in V\}.$$ Then $R$ is very clean if and
only if}\vspace{-.5mm}
\begin{enumerate}
\item [(1)]{\it $\frac{1}{2}\in A$ and $\frac{1}{2}\in B$; or} \vspace{-.5mm}
\item [(2)]{\it $R$ is strongly clean.}
\vspace{-.5mm}
\end{enumerate}\no{\it Proof.}\ \ Suppose that $R$ is very clean.
If $\frac{1}{2}\not\in A$, then $2\in J(A)$. If $a\in 1+J(A), b\in
J(B)$, then $a\pm 1\in J(A)$. By Lemma 2.1, $ax-xb=v$ is solvable.
By virtue of [10, Example 2], $R$ is strongly clean. If
$\frac{1}{2}\not\in B$, then $2\in J(B)$. If $a\in 1+J(A), b\in
J(B)$, then $a-1\in J(A)$. Further, $(b-1)\pm 1\in J(B)$. Let
$v\in V$. By virtue of Lemma 2.1, $(a-1)x-x(b-1)=v$ is solvable.
Thus, we can find some $x\in V$ such that $ax-xb=v$. In view of
[10, Example 2], $R$ is strongly clean.

We now prove the converse. If $a\pm 1\in J(A), b\in J(B)$ or $a\in
J(A), b\pm 1\in J(B)$, and $v\in V$, then $2\in J(A)$ or $2\in
J(B)$, thus $\frac{1}{2}\not\in A$ or $\frac{1}{2}\not\in B$. This
implies that $R$ is strongly clean. Therefore $R$ is very clean,
as asserted. \hfill$\Box$

\vskip4mm \no{\bf Lemma 2.3.}\ \ {\it Let $R$ be a commutative
ring with exactly two maximal ideals and suppose that
$\frac{1}{2}\in R$. Then $R$ is very clean.} \vskip2mm\no {\it
Proof.}\ \ In view of [2, Proposition 16], for any $a\in R$, there
exists an idempotent $e\in R$ such that $ea=ae$ and $a-e\in U(R)$
or $a+e\in U(R)$. Therefore $R$ is very clean.\hfill$\Box$

\vskip4mm In view of [2, Example 17], ${\Bbb Z}_{(3)}\bigcap {\Bbb
Z}_{(5)}$ is a commutative ring with exactly two maximal ideals.
We extend this result and derive the following.

\vskip4mm \no{\bf Lemma 2.4.}\ \ {\it Let $p,q\neq 2$ be prime. If $(p,q)=1$, then the ring
${\Bbb Z}_{(p)}\bigcap {\Bbb Z}_{(q)}$ is very clean, but it is not strongly clean.}
\vskip2mm\no {\it Proof.}\ \ Set $R={\Bbb Z}_{(p)}\bigcap {\Bbb Z}_{(q)}$. If $M$ is an ideal of $R$ such that
$pR\subsetneq M\subseteq R$. Choose $\frac{m}{n}\in M$ while
$\frac{m}{n}\not\in pR$. Then $p\nmid m$; hence, $(p,m)=1$. Thus, we can find some $k,l\in {\Bbb Z}$ such that $kp+lm=1$.
Clearly, $\frac{m}{1}=\frac{m}{n}\frac{n}{1}\in M$, hence, $\frac{1}{1}=p\cdot \frac{k}{1}+l\frac{m}{1}\in M$. This implies
that $M=R$. Therefore $pR$ is a maximal ideal of $R$. Likewise, $qR$ is a maximal ideal of $R$.
As $p,q\neq 2$, we see that $p,q\nmid 2$, and so $\frac{1}{2}\in R$. Obviously, $J(R)\subseteq pR\bigcap qR$. For any
$\frac{m}{n}\in pR\bigcap qR$ and $\frac{a}{b}\in R$, then
$\frac{1}{1}-\frac{m}{n}\frac{a}{b}=\frac{nb-ma}{nb}$.
Write $\frac{m}{n}=\frac{ps}{t}$. Then $psn=mt$; hence, $p~|~mt$. As $p\nmid t$, we see that $p~|~m$. Likewise,
$q~|~m$. Clearly, $p\nmid nb$, and so $p\nmid (nb-ma)$. Similarly, we see that $q\nmid (nb-ma)$.
This yields that $\frac{nb}{nb-ma}\in R$.

Thus, $\big(\frac{1}{1}-\frac{m}{n}\frac{a}{b}\big)\cdot \frac{nb}{nb-ma}=\frac{1}{1}.$
This means that $\frac{1}{1}-\frac{m}{n}\frac{a}{b}\in U(R)$, and then
$\frac{m}{n}\in J(R)$. Therefore $J(R)=pR\bigcap qR$. Assume that $M$ is a maximal ideal of $R$ and $M\neq pR, qR$.
Then $pR+M=R$ and $qR+M=R$. It follows that $R=(pR+M)(qR+M)\subseteq pR\bigcap qR+M=J(R)+M\subseteq M$, and so $R=M$. This gives a contradiction. Therefore $R$ be a commutative ring with exactly two maximal ideals. According to Lemma 2.3, $R$ is very clean.

As $(p,q)=1$, we see that $\frac{p(q+1)}{p+q}\in R$. Observing that $R$ is an integral domain, the set of
all idempotents in $R$ is $\{\frac{0}{1}, \frac{1}{1}\}$. As $q\nmid q+1$, we see that $\frac{p(q+1)}{p+q}\not\in U(R)$.
As $p\nmid (p-1)q$, we see that $\frac{p(q+1)}{p+q}-\frac{1}{1}\not\in U(R)$. This shows that $\frac{p(q+1)}{p+q}\in R$
is not strongly clean, as required.\hfill$\Box$

\vskip4mm \no{\bf Theorem 2.5.}\ \ {\it The triangular matrix ring $T_2\big(
{\Bbb Z}_{(3)}\bigcap {\Bbb Z}_{(5)}\big)$ is very clean, but it is not strongly clean.}
\vskip2mm \no{\it Proof.}\ \ As $\frac{1}{2}\in {\Bbb Z}_{(3)}\bigcap {\Bbb Z}_{(5)}\big)$, it follows by
Theorem 2.2 that $T_2\big(
{\Bbb Z}_{(3)}\bigcap {\Bbb Z}_{(5)}\big)$ is very clean, and we therefore complete the proof by Lemma 2.4.\hfill$\Box$

\vskip10mm \bc{\bf 3. $2\times 2$ Full Matrices}\ec

\vskip4mm \no The aim of this section is to investigate very cleanness of $2\times 2$ full matrices over local rings.

\vskip4mm \no{\bf Lemma 3.1.}\ \ {\it Let $R$ be a commutative
local ring, $\frac{1}{2}\in R$. Then $M_2(R)$ is very
clean.}\vskip2mm \no {\it Proof.}\ \ Let $\varphi\in M_2(R)$.
Write $\chi (\varphi)=t^2+at+b$. If $\varphi\in GL_2(R)$, or
$I_2-\varphi \in GL_2(R)$; or $I_2+\varphi \in GL_2(R)$, then
$\varphi $ is very clean. Otherwise, we may assume that
$det(\varphi), det\big(I_2-\varphi\big),
det\big(I_2+\varphi\big)\in J(R)$. Thus, $b, 1+a+b, 1-a+b\in
J(R)$. Hence, $2a\in J(R), a\in U(R)$. This implies that $2\in
J(R)$, a contradiction. Therefore we conclude that $M_2(R)$ is
very clean. \hfill$\Box$

\vskip4mm Let $R={\Bbb Z}[x]_{(x)}=\{ \frac{f(x)}{g(x)}~|~g(0)\neq
0\}$ be the localization of ${\Bbb Z}[x]$ at $(x)$. Then $R$ is a
commutative local ring with $\frac{1}{2}\in R$. It follows by
Lemma 3.1 that $M_2(R)$ is very clean.

\vskip4mm \no{\bf Lemma 3.2.}\ \ {\it Let $R$ be a ring,
$2\in J(R)$, and let $a\in R$. Then $a$ is very clean if and only
if $a$ is strongly clean.} \vskip2mm \no {\it Proof.}\ \ If $a\in
R$ is strongly clean, then it is very clean. Conversely, assume
that $a\in R$ is very clean. Then there exist an idempotent $e\in
R$ and a unit $u\in R$ such that $ae=ea$ and either $a=e+u$ or
$a=-e+u$. If $a=-e+u$, then $a=e+(u-2e)$. As $2\in J(R)$, we see
that $u-2e\in U(R)$; hence, $a\in R$ is strongly clean. Therefore
we complete the proof.\hfill$\Box$

\vskip4mm \no{\bf Theorem 3.3.}\ \ {\it Let $R$ be a commutative
local ring, and let $\varphi\in M_2(R)$. Then $\varphi\in M_2(R)$ is very clean if and only if}
\begin{enumerate}
\item [(1)]{\it $\frac{1}{2}\in R$; or}
\item [(2)]{\it $\varphi\in M_2(R)$ is strongly clean.}
\end{enumerate}
\vspace{-.5mm} \no {\it Proof.}\ Assume that $\varphi\in M_2(R)$ is very clean. If
$\frac{1}{2}\not\in R$, then $2\in J(R)$. Thus, $2I_2\in
J\big(M_2(R)\big)$. In view of Lemma 3.2, $\varphi\in
M_2(R)$ is strongly clean.

Conversely, if $\frac{1}{2}\in R$, it follows from
Lemma 3.1 that $\varphi\in M_2(R)$ is very clean. If $\varphi\in
M_2(R)$ is strongly clean, then it is very clean. Therefore $\varphi$ is
very clean in any case. \hfill$\Box$

\vskip4mm \no{\bf Example 3.4.}\ \ {\it Let $p\in {\Bbb Z}$ be a
prime, and $p\neq 2$. Then $\left(
\begin{array}{rr}
1&p\\
-1&0
\end{array}
\right)\in M_2\big({\Bbb Z}_{(p)}\big)$ is very clean, while it is
not strongly clean.}\vskip2mm \no {\it Proof.}\ \ As $p\neq 2$,
$(2,p)=1$, we can find some $k,l\in {\Bbb Z}$ such that $2k+pl=1$;
hence, $2\in U({\Bbb Z}_{(p)})$. In view of Theorem 3.3,
$M_2({\Bbb Z}_{(p)})$ is a very clean ring, and so $\left(
\begin{array}{rr}
1&p\\
-1&0
\end{array}
\right)\in M_2\big({\Bbb Z}_{(p)}\big)$ is very clean. But it is not strongly clean from
[5, Corollary 16.4.33].\hfill$\Box$

\vskip4mm For $r\in R$, define ${\Bbb S}_r=\{ f\in
R[t]~|~f~\mbox{monic, and}~f(r)\in U(R)\}.$

\vskip4mm \no{\bf Lemma 3.5.}\ \ {\it Let $R$ be a commutative
local ring, $n\geq 2$, and let $h\in R[t]$ be a monic polynomial
of degree $n$. Then the following are equivalent:}
\begin{enumerate}
\item [(1)]{\it Every $\varphi\in M_n(R)$ with $\chi (\varphi)=h$ is very clean.}
\item [(2)]{\it There exists a factorization $h=h_0h_1$ such that
$(h_0,h_1)=R[t]$ and $h_0\in {\Bbb S}_0$ and $h_1\in {\Bbb
S}_1\bigcup {\Bbb S}_{-1}$.}
\end{enumerate}
\vspace{-.5mm} \no {\it Proof.}\ \  $(1)\Rightarrow (2)$ Since
$\varphi$ is very clean, we see that $\varphi$ or $-\varphi$ is
strongly clean. If $\varphi$ is strongly clean, it follows by [4,
Theorem 12] that there exists a factorization $h=h_0h_1$ such that
$(h_0,h_1)=R[t]$ and $h_0\in {\Bbb S}_0$ and $h_1\in {\Bbb S}_1$.

If $-\varphi$ is strongly clean, it follows by [4, Theorem 12]
that $g(\mu):=det\big(\mu I_n-(-\varphi)\big)=g_0g_1$ where
$(g_0,g_1)=R[\mu]$ and $g_0\in {\Bbb S}_0$ and $g_1\in {\Bbb
S}_1$. This implies that $$h(t)=det\big(t
I_n-\varphi\big)=(-1)^ng(-t)=(-1)^{n}g_0(-t)g_1(-t).$$ Set
$h_0=(-1)^{degg_0}g_0(-t)$ and $h_1=(-1)^{degg_1}g_1(-t)$. Then
$h=h_0h_1$ with $(h_0,h_1)=R[t]$. Further,
$h_0(0)=(-1)^{degg_0}g_0(0)\in U(R)$; hence, $h_0\in {\Bbb S}_0$.
In addition, $h_1(-1)=(-1)^{degg_1}g_1(1)\in U(R)$. This implies
that $h_1\in {\Bbb S}_{-1}$. Therefore  $h_0\in {\Bbb S}_0$ and
$h_1\in {\Bbb S}_1\bigcup {\Bbb S}_{-1}$.

$(2)\Rightarrow (1)$ By hypothesis, there exists a factorization
$h=h_0h_1$ such that $(h_0,h_1)=R[t]$ and $h_0\in {\Bbb S}_0$ and
$h_1\in {\Bbb S}_1\bigcup {\Bbb S}_{-1}$. If $h_1\in {\Bbb S}_1$, it
follows by [4, Theorem 12] that $\varphi\in M_n(R)$ is strongly
clean, and so it is very clean. If $h_1\in {\Bbb S}_{-1}$, then
$h_1(-1)\in U(R)$. Set $g(\mu):=(-1)^nh(-\mu)$. Then
$g(\mu):=g_0g_1$ where $g_0(\mu)=(-1)^{degh_0}h_0(-\mu)$ and
$g_1(\mu)=(-1)^{degh_1}h_1(-\mu)$. As $g_0(0)\in U(R)$, we see
that $g_0\in {\Bbb S}_0$. Further, $g_1(1)=(-1)^{degh_1}h_1(-1)\in
U(R)$, and then $g_1\in {\Bbb S}_1$. Clearly, $g(\mu)=det\big(\mu
I_n-(-\varphi)\big)$. In view of [4, Theorem 12], $-\varphi\in
M_n(R)$ is strongly clean. Therefore $\varphi\in M_n(R)$ is very
clean, as asserted.\hfill$\Box$

\vskip4mm In what follows, we consider more explicit criteria for
very clean $2\times 2$ matrices over commutative rings.

\vskip4mm \no{\bf Theorem 3.6.}\ \ {\it Let $R$ be a commutative
local ring, and let $h\in R[t]$ be a monic polynomial of degree
$2$. Then the following are equivalent:}
\begin{enumerate}
\item [(1)]{\it Every $\varphi\in M_2(R)$ with $\chi (\varphi)=h$ is very clean.}
\item [(2)]{\it There exists a factorization $h=h_0h_1$ such that $h_0\in {\Bbb S}_0$ and $h_1\in {\Bbb
S}_1\bigcup {\Bbb S}_{-1}$.}
\end{enumerate}
\vspace{-.5mm} \no {\it Proof.}\ $(1)\Rightarrow (2)$ is trivial
from Lemma 3.5.

$(2)\Rightarrow (1)$ By hypothesis, there exists a factorization
$\chi(\varphi )=h_0h_1$ such that $h_0\in {\Bbb S}_0$ and $h_1\in
{\Bbb S}_1\bigcup {\Bbb S}_{-1}$.

Case I. $deg(h_0)=2$ and $deg(h_1)=0$. Then $h_0=
t^2-tr(\varphi)t+det(\varphi)$ and $h_1=1$. Hence, $(h_0,h_1)=R[t]$.

Case II. $deg(h_0)=1$ and $deg(h_1)=1$. Then $h_0=t-\alpha$ and
$h_1=t-\beta$. Since $h_0\in {\Bbb S}_0$, $\alpha\in
U(R)$. As $h_1\in {\Bbb S}_1\bigcup {\Bbb S}_{-1}$, we see that
$\beta\in 1+U(R)$ or $\beta\in -1+U(R)$.

If $\beta-\alpha\in U(R)$, then
$h_0(\beta-\alpha)^{-1}-h_1(\beta-\alpha)^{-1}=1$; hence, $(h_0,h_1)=R[t]$.

If $\beta-\alpha\in J(R)$, then $\beta\in U(R)$. Let $h_0'=t^2-(\alpha+\beta)t+\alpha\beta$ and $h_1'=1$. Then
$h_0'\in {\Bbb S}_0$ and $h_1'\in {\Bbb S}_1$. In addition, $(h_0',h_1')=R[t]$.

Case III. $deg(h_0)=0$ and $deg(h_1)=2$. Then
$h_0=1$ and $h_1=det\big(tI_2-\varphi\big)$. Hence, $(h_0,h_1)=R[t]$.

In any case, there exists a factorization $h=h_0h_1$ such that $h_0\in {\Bbb S}_0, h_1\in {\Bbb
S}_1\bigcup {\Bbb S}_{-1}$ and $(h_0,h_1)=R[t]$. Therefore we complete the proof by Lemma 3.5.\hfill$\Box$

\vskip4mm \no{\bf Corollary 3.7.}\ \ {\it Let $R$ be a commutative
local ring, and let $\varphi\in M_2(R)$. Then $\varphi$ is very
clean if and only if $\varphi\in M_2(R)$ is strongly clean; or
$I_2+\varphi \in GL_2(R)$.}\vskip2mm\no {\it Proof.}\ \ Assume
that $\varphi\in M_2(R)$ is strongly clean; or $I_2+\varphi \in
GL_2(R)$. Then $\varphi$ is very clean. Conversely, assume that
$\varphi$ is very clean and $I_2+\varphi \not\in GL_2(R)$. We may
assume that $\varphi, I_2-\varphi\not\in GL_2(R)$, and so
$det(\varphi), det\big(I_2-\varphi\big),
det\big(I_2+\varphi\big)\in J(R)$. It follows from
$det(\varphi)=-\alpha(\beta +a)\in J(R)$ that $\beta+a\in J(R)$.
Hence, $1+\beta +a\in U(R)$. Set $h_0=t-\alpha$ and $h_1=t+\beta
+a$. Then $h_0\in {\Bbb S}_0, h_1\in {\Bbb S}_1$. As in the proof
of Theorem 3.6, we may assume that $(h_0,h_1)=1$. In light of [4,
Theorem 12], $\varphi$ is clean, as required.\hfill$\Box$

\vskip15mm \bc{\bf 4. MATRICES OVER POWER SERIES}\ec

\vskip4mm \no Let $a\in R$. Then $l_a: R\to R$ and $r_a: R\to R$
denote, respectively, the abelian group endomorphisms given by
$l_a(r)=ar$ and $r_a(r)=ra$ for all $r\in R$. Thus, $l_a-r_b$ is
an abelian group endomorphism such that $(l_a-r_b)(r)=ar-rb$ for
any $r\in R$. Following Diesl, a local ring $R$ is weakly bleached
provided that for any $a\in 1+J(R), b\in J(R)$, ${\ell}_a-r_b,
{\ell}_b-r_a: R\to R$ are surjective. The class of weakly bleached
local rings contains many familiar examples, e.g., commutative
local rings, local rings with nil Jacobson radicals, local rings
for which some power of each element of their Jacobson radicals is
central (cf. [3, Example 13]). The goal of this section is to
investigate very clean matrices with power series entries over
local rings.

\vskip4mm \no{\bf Lemma 4.1.}\ \ {\it Let $R$ be a ring. Then
$f(x)\in R[[x]]$ is very clean if there exists an idempotent $e\in
R$ such that} \vspace{-.5mm}
\begin{enumerate}
\item [(1)]{\it $f(0)e=ef(0)$;}
\vspace{-.5mm}
\item [(2)]{\it $f(0)-e\in U(R)$ or $f(0)+e\in U(R)$; }\vspace{-.5mm}
\item [(3)]{\it for any $b\in R$, there exists an $x\in R$ such that $[f(0),x]=[e,b]$.}
\end{enumerate}\no {\it Proof.}\ \ By hypothesis, there exists an
idempotent $e\in R$ such that $f(0)e=ef(0)$ and either $f(0)-e\in
U(R)$ or $f(0)+e\in U(R)$. Assume that $f(0)-e\in U(R)$. For any
$b\in R$, we can find an $x\in R$ such that $[f(0),x]=[e,b]$.
Applying [8, Theorem 3.3.2] to $f(x)$, we see that $f(x)\in
R[[x]]$ is strongly clean. If $f(0)+e\in U(R)$, then $-f(0)-e\in
U(R)$. For any $b\in R$, we can find an $x\in R$ such that
$[f(0),x]=[e,b]$. Choose $y=-x$. Then $[-f(0),y]=[e,b]$. Applying
[8, Theorem 3.3.2] to $-f(x)$, we see that $-f(x)\in R[[x]]$ is
strongly clean. Therefore $f(x)$ is very clean, as
required.\hfill$\Box$

\vskip4mm \no{\bf Theorem 4.2.}\ \ {\it Let $R$ be a weakly
bleached local ring. Then the following are equivalent:}
\begin{enumerate}
\item [(1)]{\it $A(x)\in M_2\big(R[[x]])$ is very clean.}
\vspace{-.5mm}
\item [(2)]{\it $A(0)\in M_2(R)$ is very clean.}\vspace{-.5mm}
\end{enumerate}\no {\it Proof.}\ \ $(1)\Rightarrow (2)$ Since $A(x)$ is very clean in
$M_2\big(R[[x]]\big)$, there exist an $E(x)=E^2(x)\in
M_2\big(R[[x]]\big)$ and a $U(x)\in GL_2(R[[x]]\big)$ such that
$E(x)U(x)=U(x)E(x)$, and that either $A(x)=E(x)+U(x)$ or
$A(x)=-E(x)+U(x)$. This implies that $E(0)U(0)=U(0)E(0)$, and that
either $A(0)=E(0)+U(0)$ or $A(0)=-E(0)+ U(0)$, where
$E(0)=E^2(0)\in M_2(R)$ and $U(0)\in GL_2(R)$. As a result, $A(0)$
is very clean in $M_2(R)$.

$(2)\Rightarrow (1)$ Clearly, $M_2\big(R[[x]]\big)\cong
M_2(R)[[x]]$. We may assume that $A(x)\in M_2(R)[[x]]$. Set
$S=M_2(R)$. Then $A(0)= \left(
\begin{array}{cc}
a&b\\
c&d
\end{array}
\right)\in S$. Then there exists an idempotent $e=\left(
\begin{array}{cc}
e_{11}&e_{12}\\
e_{21}&e_{22}
\end{array}
\right)$ $\in S$ such that $A(0)e=eA(0)$ and either $A(0)-e\in
U(S)$ or $A(0)+e\in U(S)$. In view of [5, Lemma 16.4.10], there
exists some $u\in U(S)$ such that $ueu^{-1}=\left(
\begin{array}{cc}
e_1&0\\
0&e_2
\end{array}
\right)$. Clearly, $e_1=e_1^2, e_2=e_2^2$. As $R$ is local, $e_1$
and $e_2$ are trivial idempotents. If $e_1=e_2=0$ or $e_1=e_2=1$,
then $e=0$ or $e=I_2$, and so for any $s\in S$, there exists an
$x=0$ such that $[A(0),x]=[e,s]$. Thus, we may assume that $e_1=1$
and $e_2=0$. It follows from $A(0)e=eA(0)$ that
$\big(uA(0)u^{-1}\big)(ueu^{-1})=(ueu^{-1})\big(uA(0)u^{-1}\big)$;
hence, $uA(0)u^{-1}=\left(
\begin{array}{cc}
a_{11}&0\\
0&a_{22}
\end{array}
\right)$, where $a_{22}\in U(R)$ and either $a_{11}\in 1+U(R)$ or
$a_{11}\in -1+U(R)$. Set $\alpha:=\left(
\begin{array}{cc}
a_{11}&0\\
0&a_{22}
\end{array}
\right)$. Assume that $a_{11}\in 1+U(R)$ and $a_{22}\in U(R)$. If
$a_{22}\in 1+U(R)$, then we choose $f=I_2$, then $\alpha -f\in
U(S), \alpha f=f\alpha$ and that for any $\beta\in S$,
$[\alpha,0]=[f,\beta]$. If $a_{11}\in U(R)$, then we choose $f=0$,
then $\alpha -f\in U(S), \alpha f=f\alpha$ and that for any
$\beta\in S$, $[\alpha,0]=[f,\beta]$. Thus, we assume that
$a_{11}\in J(R), a_{22}\in 1+J(R)$. Choose $f=\left(
\begin{array}{cc}
1&0\\
0&0
\end{array}
\right)\in R$. Then $\alpha -f\in U(S), \alpha f=f\alpha$. For any
$\beta=(\beta_{ij})\in S$, as $R$ is weakly bleached, there exist
some $x_1,x_2\in S$ such that $a_{11}x_1-x_1a_{22}=\beta_{12}$ and
$a_{22}x_2-x_2a_{11}=-\beta_{21}$. Choose $x=\left(
\begin{array}{cc}
0&x_1\\
x_2&0
\end{array}
\right)\in R$. It is easy to verify that $$ [\alpha,x]=\left(
\begin{array}{cc}
0&a_{11}x_1-x_1a_{22}\\
a_{22}x_2-x_2a_{11}&0
\end{array}
\right)=\left(
\begin{array}{cc}
0&\beta_{12}\\
-\beta_{21}&0
\end{array}
\right)=[f,\beta].$$

Assume that $a_{11}\in -1+U(R)$ and $a_{22}\in U(R)$. If
$a_{22}\in -1+U(R)$, then we choose $f=I_2$, then $\alpha +f\in
U(S), \alpha f=f\alpha$ and that for any $\beta\in S$,
$[\alpha,0]=[f,\beta]$. If $a_{11}\in U(R)$, then we choose $f=0$,
then $\alpha -f\in U(S), \alpha f=f\alpha$ and that for any
$\beta\in S$, $[\alpha,0]=[f,\beta]$. Thus, we assume that
$a_{11}\in J(R), a_{22}\in -1+J(R)$. Thus, $-a_{11}\in J(R),
-a_{22}\in 1+J(R)$. Choose $f=\left(
\begin{array}{cc}
1&0\\
0&0
\end{array}
\right)\in S$. Then $\alpha +f\in U(S), \alpha f=f\alpha$. For any
$\beta=(\beta_{ij})\in S$, as $R$ is weakly bleached, there exist
some $x_1,x_2\in R$ such that
$(-a_{11})x_1-x_1(-a_{22})=-\beta_{12}$ and
$(-a_{22})x_2-x_2(-a_{11})=\beta_{21}$. Choose $x=\left(
\begin{array}{cc}
0&x_1\\
x_2&0
\end{array}
\right)\in S$. Then we check that $ [\alpha,x]=\left(
\begin{array}{cc}
0&a_{11}x_1-x_1a_{22}\\
a_{22}x_2-x_2a_{11}&0
\end{array}
\right)=\left(
\begin{array}{cc}
0&\beta_{12}\\
-\beta_{21}&0
\end{array}
\right)=[f,\beta]. $ The case $e_1=0$, $e_2=1$ is similar.

For any $s\in S$, it follows from the preceding discussion that
there exists an $x\in S$ such that
$[uA(0)u^{-1},x]=[ueu^{-1},usu^{-1}]$. Therefore $[A(0),
u^{-1}xu]=[e,s]$. According to Lemma 4.1, $A(x)$ is very
clean.\hfill$\Box$

\vskip4mm Let $p~ (\neq 2)$ be a prime number. In light of Theorem
4.2, the ring $M_2\big({\Bbb Z}_{(p)}[[x]]\big)$ is very clean.
But $M_2\big({\Bbb
Z}_{(p)}[[x]]\big)$ is not strongly clean.

\vskip4mm \no{\bf Corollary 4.3.}\ \ {\it Let $R$ be a commutative
local ring, and let $A(x)\in M_2\big(R[[x]])$. Then the following
are equivalent:} \vspace{-.5mm}
\begin{enumerate}
\item [(1)]{\it $A(x)\in M_2\big(R[[x]])$ is very clean.}
\vspace{-.5mm}
\item [(2)]{\it $A(0)\in M_2(R)$ is very clean.}\vspace{-.5mm}
\end{enumerate}\no {\it Proof.}\ \ Since every commutative local ring is weakly bleached, we complete the proof by
Theorem 4.2.\hfill$\Box$

\vskip4mm \no{\bf Corollary 4.4.}\ \ {\it Let $R$ be a commutative
local ring, and let $A(x)\in M_2\big(R[[x]]/(x^{m})\big)$ $(m\geq
1)$. Then the following are equivalent:} \vspace{-.5mm}
\begin{enumerate}
\item [(1)]{\it $A(x)\in M_2\big(R[[x]]/(x^m)\big)$ is very clean.}
\vspace{-.5mm}
\item [(2)]{\it $A(0)\in M_2(R)$ is very clean.}\vspace{-.5mm}
\end{enumerate}\no {\it Proof.}\ \ It is obvious from Corollary 4.3.\hfill$\Box$

\vskip4mm \no{\bf Example 4.5.}\ \ {\it Let $R={\Bbb
Z}_4[x]/(x^2)$, and let $$A(x)=\left( \begin{array}{cc}
\overline{3}&\overline{2}+\overline{2}x\\
\overline{2}+x&\overline{3}x \end{array} \right)\in M_2(R).$$
Obviously, ${\Bbb Z}_4$ is a commutative local ring, and that
$R={\Bbb Z}_4[[x]]/(x^2)$. Since we have the very clean
decomposition $A(0)=\left(
\begin{array}{cc}
\overline{3}&\overline{0}\\
\overline{0}&\overline{3} \end{array} \right)+\left(
\begin{array}{cc}
\overline{0}&\overline{2}\\
\overline{2}&\overline{1} \end{array} \right)$ in $M_2({\Bbb
Z}_4)$, it follows by Corollary 4.4 that $A(x)\in M_2(R)$ is very
clean.}\hfill$\Box$

\vskip4mm \no{\bf Theorem 4.6.}\ \ {\it Let $R$ be a weakly
bleached local ring. Then the following are equivalent:}
\vspace{-.5mm}
\begin{enumerate}
\item [(1)]{\it $A(x)\in T_2\big(R[[x]])$ is very clean.}
\vspace{-.5mm}
\item [(2)]{\it $A(0)\in T_2(R)$ is very clean.}\vspace{-.5mm}
\end{enumerate}\no {\it Proof.}\ \ $(1)\Rightarrow (2)$ is obvious.

$(2)\Rightarrow (1)$ Let $A(0)= \left(
\begin{array}{cc}
a_1&c\\
0&a_2
\end{array}
\right)\in R$ and $S=T_2(R)$.\\ $(a)$ $a_1, a_2\in 1+U(R)$. Then we
choose $e=\left(
\begin{array}{cc}
1&0\\
0&1
\end{array}
\right)$. Then $A(0)-e\in U(S)$ and $A(0)e=eA(0)$. For any $b\in
S$, we choose $x=0$. Then $[A(0),x]=[e,b]$.

$(b)$ $a_1, a_2\in -1+U(R)$. Then we choose $e=\left(
\begin{array}{cc}
1&0\\
0&1
\end{array}
\right)$. Then $A(0)+e\in U(S)$ and $A(0)e=eA(0)$. For any $b\in
S$, we choose $x=0$. Then $[A(0),x]=[e,b]$.

$(c)$ $a_1, a_2\in U(R)$. Then $A(0)\in U(S)$, and so we choose
$e=0$. Then $A(0)-e\in U(S)$ and $A(0)e=eA(0)$. For any $b\in S$,
we choose $x=0$. Then $[A(0),x]=[0,b]$.

Thus, either $a_1\in 1+J(R), a_2\in J(R)$ or $a_1\in J(R), a_2\in
1+J(R)$, and that either $a_1\in -1+J(R), a_2\in J(R)$ or $a_1\in
J(R), a_2\in -1+J(R)$. Therefore we may assume that either $a_1\in
\pm 1+J(R), a_2\in J(R)$ or $a_1\in J(R), a_2\in \pm 1+J(R)$. For
such $a\in S$, by hypothesis, there exist an idempotent $e\in S$
and a unit $u\in S$ such that $A(0)e=eA(0)$ and either $a-e=u$ or
$a+e=u$.

Case I. $a_1\in \pm 1+J(R),a_2\in J(R)$. Then $e=\left(
\begin{array}{cc}
0&y\\
0&1
\end{array}
\right)$. For any $b=\left(
\begin{array}{cc}
b_1&z\\
0&b_2
\end{array}
\right)\in S$, as $R$ is weakly bleached, we can find a $w\in R$
such that $a_1w-wa_2=yb_2-b_1y-z$. Choose $x=\left(
\begin{array}{cc}
0&w\\
0&0
\end{array}
\right)\in S$. Then $[A(0),x]=\left(
\begin{array}{cc}
0&a_1w-wa_2\\
0&0
\end{array}
\right)$ $=\left(
\begin{array}{cc}
0&yb_2-b_1y-z\\
0&b_2
\end{array}
\right)=[e,b]$.

Case II. $a_1\in J(R),a_2\in \pm 1+J(R)$. Then $e=\left(
\begin{array}{cc}
1&y\\
0&0
\end{array}
\right)$. For any $b=\left(
\begin{array}{cc}
b_1&z\\
0&b_2
\end{array}
\right)\in S$, as $R$ is weakly bleached, we can find a $w\in R$
such that $a_1w-wa_2=z+yb_2-b_1y$. Choose $x=\left(
\begin{array}{cc}
0&w\\
0&0
\end{array}
\right)\in S$. Then $$[A(0),x]=\left(
\begin{array}{cc}
0&a_1w-wa_2\\
0&0
\end{array}
\right)=\left(
\begin{array}{cc}
0&z+yb_2-b_1y\\
0&0
\end{array}
\right)=[e,b].$$ Thus, $A(0)$ satisfies the conditions in Lemma
4.1. Since $T_2\big(R[[x]]\big)\cong T_2(R)[[x]]$, we conclude
that $A(x)$ is very clean.\hfill$\Box$

\vskip4mm \no{\bf Corollary 4.7.}\ \ {\it Let $R$ be a local ring.
Then the following are equivalent:}
\begin{enumerate}
\item [(1)]{\it $T_2(R)$ is very clean.}
\vspace{-.5mm}
\item [(2)]{\it $T_2\big(R[[x]]\big)$ is very clean.}
\end{enumerate}\no
{\it Proof.}\ \ $(1)\Rightarrow (2)$ Since $T_2(R)$ is very clean, it follows from Theorem 2.2
that $\frac{1}{2}\in R$ or $T_2(R)$ is strongly clean. If
$\frac{1}{2}\in R$, as in the proof of Theorem 2.2, $T_2(R)$ is
very clean. If $T_2(R)$ is strongly clean, it follows from [10,
Example 2] that $R$ is weakly bleached. Therefore
$T_2\big(R[[x]]\big)$ is very clean by Theorem 4.6.

$(2)\Rightarrow (1)$ is just as easy.\hfill$\Box$

\vskip4mm Let $R$ be a commutative local ring. If $\frac{1}{2}\in
R$, then $T_2\big(R[[x]]\big)$ is very clean. As in the proof of
Theorem 2.2, $T_2(R)$ is very clean. Therefore we are done by
Corollary 4.7. For instance, for any prime $p~ (\neq 2)$. Then
$T_2\big({\Bbb Z}_{(p)}[[x]]\big)$ is very clean.

\vskip10mm \bc{\large\bf REFERENCES}\ec \vskip 4mm \small{\re{1}
M.S. Ahn and D.D. Anderson, Weakly clean rings and almost clean
rings, {\it Rocky Mountain J. Math.}, {\bf 36}(2006), 783--798.

\re{2} D.D. Anderson and V.P. Camillo, Commutative rings whose
elements are a sum of a unit and idempotent, {\it Comm. Algebra},
{\bf 30}(2002), 3327--3336.

\re{3} G. Borooah; A.J. Diesl and T.J. Dorsey, Strongly clean
triangular matrix rings over local rings, {\it J. Algebra}, {\bf
312}(2007), 773--797.

\re{4} G. Borooah; A.J. Diesl and T.J. Dorsey, Strongly clean
matrix rings over commutative local rings, {\it J. Pure Appl.
Algebra}, {\bf 212}(2008), 281--296.

\re{5} H. Chen, {\it Rings Related Stable Range Conditions},
Series in Algebra 11, World Scientific, Hackensack, NJ, 2011.

\re{6} H. Chen, On strongly nil clean matrices, {\it Comm.
Algebra}, {\bf 41}(2013), 1074-1086.

\re{7} T.J. Dorsey, {\it Cleanness and Strong Cleanness of Rings
of Matrices}, Ph.D. Thesis, University of California, Berkeley,
2006.

\re{8} D.R. Shifflet, {\it Optimally Clean Rings}, Ph.D. Thesis,
Bowling Green State University, 2011.

\re{9} W.K. Nicholson, Lifting idempotents and exchange rings,
{\it Trans. Amer. Math. Soc.}, {\bf 229}(1977), 269--278.

\re{10} W.K. Nicholson, Strongly clean rings and Fitting's lemma,
{\it Comm. Algebra}, {\bf 27}(1999), 3583--3592.

\re{11} W.K. Nicholson, Clean rings: a survey, {\it Advances in
Ring Theory}, World Sci. Publ., Hackensack, NJ, 2005, 181--198.

\end{document}